\documentclass[letterpaper, 10 pt, conference]{ieeeconf}
\IEEEoverridecommandlockouts

\usepackage[utf8]{inputenc}
\usepackage{amsmath}
\usepackage{mathrsfs}
\usepackage{amssymb}
\usepackage{color}
\usepackage{dsfont}
\usepackage{cite}
\usepackage{comment}
\usepackage{float}
\usepackage{graphicx}
\usepackage[font = footnotesize]{caption}
\usepackage[font = footnotesize]{subcaption}
\usepackage{bernsteinStyle}
\usepackage{float}
\usepackage{tikz}
\usepackage{booktabs}
\usepackage{tabularx}
\usepackage{multirow}
\usepackage{bm}
\usepackage{makecell}
\usepackage{indentfirst}
\newcolumntype{C}[1]{>{\centering\let\newline\\\arraybackslash\hspace{0pt}}m{#1}}

\usepackage{color, colortbl}
\definecolor{Yellow}{rgb}{1,1,0}

\usepackage{hhline}

\usetikzlibrary{shapes,arrows,calc,positioning}
\tikzstyle{bigblock} = [draw, fill=blue!20, rectangle, 
    minimum height=6em, minimum width=8em]
\tikzstyle{medblock} = [draw, fill=blue!20, rectangle, 
    minimum height=4em, minimum width=4em]    
\tikzstyle{mux} = [draw, fill=black!20, rectangle, 
    minimum height=5em, minimum width=0.1em]    
\tikzstyle{smallblock} = [draw, fill=blue!20, rectangle, 
    minimum height=3em, minimum width=4em]
\tikzstyle{sum} = [draw, fill=blue!20, circle, node distance=1cm]
\tikzstyle{signal} = [coordinate]
\tikzstyle{pinstyle} = [pin edge={to-,thin,black}]
\tikzstyle{block} = [draw, fill=blue!20, rectangle, 
    minimum height=3em, minimum width=6em]
\tikzstyle{blockS} = [draw, fill=blue!20, rectangle, 
    minimum height=3em, minimum width=4em]  
\tikzstyle{sum} = [draw, fill=blue!20, circle, node distance=1.5cm]
\tikzstyle{gain} = [draw, fill=blue!20, regular polygon, regular polygon sides = 3, node distance=1.25cm, shape border rotate = -90]
\tikzstyle{mult} = [draw, fill=blue!20, circle, node distance=1.25cm ,inner sep=0pt, minimum size = 0.3cm]

\tikzstyle{input} = [coordinate]
\tikzstyle{output} = [coordinate]

\newcounter{example}

\title{\LARGE Retrospective Cost-based Extremum Seeking Control with Vanishing Perturbation for Online Output Minimization}

\author{\large Juan A. Paredes,  Jhon Manuel Portella, Dennis S. Bernstein, and Ankit Goel
\thanks{Juan A. Paredes and Dennis S. Bernstein are with the Department of Aerospace Engineering, University of Michigan, Ann Arbor, MI, USA. {\tt\small \{jparedes, dsbaero\}@umich.edu}}
\thanks{Jhon Manuel Portella and Ankit Goel are with the Department of Mechanical Engineering, University of Maryland, Baltimore County, MD 21250.
{\tt \small \{jportella,ankgoel\}@umbc.edu}}
}

\begin{document}

\maketitle

\begin{abstract}
Extremum seeking control (ESC) constitutes a powerful technique for online optimization with theoretical guarantees for convergence to the neighborhood of the optimizer under well-understood conditions.
However, ESC requires a nonconstant perturbation signal to provide persistent excitation to the target system to yield convergent results, which usually results in steady state oscillations. 
While certain techniques have been proposed to eliminate perturbations once the neighborhood of the minimizer is reached, system modifications and environmental perturbations can suddenly change the minimizer and nonconstant perturbations would once more be required to convergence to the new minimizer.
Hence, this paper develops a retrospective cost-based ESC (RC/ESC) technique for online output minimization with a vanishing perturbation, that is, a perturbation that becomes zero as time increases independently from the state of the controller or the controlled system.
The performance of the proposed algorithm is illustrated via numerical examples.
\end{abstract}

\section{Introduction}\label{sec:introduction}

Extremum seeking control (ESC) is a powerful adaptive control technique that leverages persistent system excitation to search for extrema in order to either minimize or maximize a used-defined metric \cite{KrsticBookESC}.
The stability and convergence properties of ESC and their conditions have been thoroughly studied and are well understood \cite{krstic2000,krstic2000_2,tan2010extremum}.
ESC has been applied in a wide arrange of fields, including robotics \cite{matveev2011,calli2012,matveev2014,vweza2015,matveev2015,bagheri2018,calli2018,sotiropoulos2019} and energy management \cite{creaby2009,johnson2012,ghaffari2014,ye2016,bizon2017,zhou2017,zhou2018}.

A feature of ESC is a persistent perturbation signal, which enables gradient estimation algorithms to yield a search direction that points towards local extrema, thus enabling convergence. 
However, the implementation of this perturbation signal results in steady state oscillations, which may be prohibitive in physical testing.
Modifications to the ESC algorithm have been proposed to address this issue, which include modifying the perturbation signal to vanish depending on controller and system values and implementing dynamics that suppress the perturbation signal once a neighborhood of the minimizer has been reached \cite{wang2016,haring2016,haring2018,yin2019,bhattacharjee2021}.

The contribution of this paper is thus an ESC algorithm for online output minimization with a vanishing perturbation, that is, a perturbation that becomes zero as time increases independently from the state of the controller or the controlled system.
Hence, the perturbation is independent from the rest of the system and no extra dynamics are required to suppress the perturbation.
In particular, we consider retrospective cost adaptive control (RCAC), which re-optimizes the coefficients of the feedback controller at each step \cite{rahmanCSM2017}.
A similar retrospective cost algorithm was proposed in \cite{goel2019,goel2023}, in which a fixed target model was used to issue a search direction in the input space for optimization.
In this work, a Kalman Filter (KF) is used to estimate the gradient of the system output which is encoded into the target model to provide a time-varying search direction to RCAC.
Hence, the combination of the KF, the target model construction procedure, and RCAC yields retrospective cost based ESC (RC/ESC).
A preliminary version of this algorithm was considered in \cite{juanScitech2022}, in which gradient estimation was performed by using a simple high-pass filter.
 
The contents of the paper are as follows.
Section \ref{sec:control_statement} provides a statement of the control problem, which involves continuous-time dynamics under sampled-data feedback control.
Section \ref{sec:ESC} provides a review of continuous-time ESC.
Section \ref{sec:RCAC_ESC} introduces RC/ESC, in which a KF estimates the gradient of the system output and provides a target model to RCAC.
Section \ref{sec:numerical_examples} presents numerical examples that illustrate the performance of RC/ESC, including examples with static maps and a example with a dynamic system.
Finally, \ref{sec:conclusions} presents conclusions.

{\bf Notation:}
$\bfq\in\BBC$ denotes the forward-shift operator.
$x_{(i)}$ denotes the $i$th component of $x\in\BBR^n.$


\section{Problem Statement}\label{sec:control_statement}

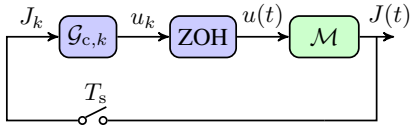
\begin{figure} [h!]
    \centering
    \resizebox{0.65\columnwidth}{!}{%
    \begin{tikzpicture}[>={stealth'}, line width = 0.25mm]

    \node [input, name=ref]{};
    \node [smallblock, rounded corners, right = 0.5cm of ref , minimum height = 0.6cm, minimum width = 0.7cm] (controller) {$\SG_{\rmc, k}$};
    \node [smallblock, rounded corners, right = 0.75cm of controller, minimum height = 0.6cm , minimum width = 0.5cm] (DA) {ZOH};
    
    \node [smallblock, fill=green!20, rounded corners, right = 0.75cm of DA, minimum height = 0.6cm , minimum width = 1cm] (system) {$\SM$};
    \node [output, right = 0.5cm of system] (output) {};
    \node [input, below = 0.9cm of system] (midpoint) {};
    
    \draw [->] (controller) -- node [above] {$u_k$} (DA);\
    \draw [->] (DA) -- node [above] {$u (t)$} (system);
    
    \node[circle,draw=black, fill=white, inner sep=0pt,minimum size=3pt] (rc11) at ([xshift=-3.2cm]midpoint) {};
    \node[circle,draw=black, fill=white, inner sep=0pt,minimum size=3pt] (rc21) at ([xshift=-3.5cm]midpoint) {};
    \draw [-] (rc21.north east) --node[below,yshift=.55cm]{$T_\rms$} ([xshift=.3cm,yshift=.15cm]rc21.north east) {};
    
    \draw [->] (system) -- node [name=y, near end]{} node [very near end, above] {$J (t)$}(output);
    
    \draw [-] (y.west) |- (midpoint);
    \draw [-] (midpoint) -| (rc11.east);
    \draw [->] (rc21) -| ([xshift = -0.75cm]controller.west) -- node [near end, above, xshift = -0.2cm] {$J_k$} (controller.west);
    
    \end{tikzpicture}
    }  
    \caption{Sampled-data implementation of adaptive controller for control of the continuous-time system $\SM.$
    All sample-and-hold operations are synchronous.
    The adaptive controller uses $J_k$ as an input and generates the discrete-time control $u_k$ at each step $k$.
    Note that all components of $J_k$ are nonnegative.
    The resulting continuous-time control $u(t)$ is generated by applying a zero-order-hold operation to $u_k$.
    The objective of the controller is yield an input signal that minimizes the output of the continuous-time system, that is, yield $u(t)$ such that $\lim_{t \to \infty} J(t) = 0.$
    }
    \label{AC_CT_blk_diag}
\end{figure}

We consider continuous-time dynamics under sampled-data control using discrete-time adaptive controllers to reflect the practical implementation of digital controllers for physical systems.
In particular, we consider the control architecture shown in Figure \ref{AC_CT_blk_diag}, where $\SM$ is the target continuous-time system, $t\ge 0$, $u(t)\in\BBR^m$ is the control, and $J(t)\in\BBR^p$ is the output of $\SM,$ whose components are all non-negative.

The output $J(t)$ is sampled to generate the measurement $J_k \in \BBR^p,$ which, for all $k\ge0,$ is given by
\begin{equation}
    J_k \isdef  J(k T_\rms),
\end{equation}
where  $T_\rms>0$ is the sample time.
The adaptive controller, which is updated at each step $k,$ is denoted by $\SG_{\rmc,k}$.
The input to $\SG_{\rmc,k}$ is $J_k$, and its output at each step $k$ is the discrete-time control $u_k\in\BBR^m.$
The continuous-time control signal $u(t)$ applied to the structure is generated by applying a zero-order-hold operation to $u_k,$ that is,
for all $k\ge0,$ and, for all $t\in[kT_\rms, (k+1) T_\rms),$ 
\begin{equation}
    u(t) = u_k.
\end{equation}
The objective of the adaptive controller is to yield an input signal that minimizes the output of the continuous-time system, that is, yield $u(t)$ such that $\lim_{t \to \infty} J(t) = 0.$

\section{Overview of Extremum Seeking Control}\label{sec:ESC}

\subsection{General Scheme}

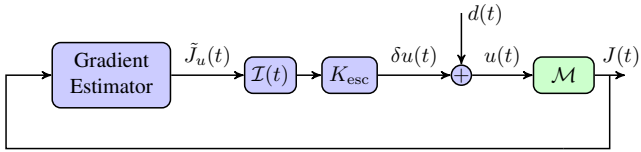
\begin{figure} [h!]
    \centering
    \resizebox{\columnwidth}{!}{%
    \begin{tikzpicture}[>={stealth'}, line width = 0.25mm]

    \node [input, name=ref]{};
    \node [smallblock, rounded corners, right = 0.5cm of ref , minimum height = 0.6cm, minimum width = 0.7cm] (estimator) {$\begin{array}{c} \mbox{Gradient} \\ \mbox{Estimator} \end{array}$};
    \node [smallblock, rounded corners, right = 1.2cm of estimator , minimum height = 0.6cm, minimum width = 0.7cm] (integrator) {$\SI (t)$};
    \node [smallblock, rounded corners, right = 0.4cm of integrator , minimum height = 0.6cm, minimum width = 0.7cm] (int_gain) {$K_{\rm esc}$};
    \node [sum, right =1.2cm of int_gain] (sum1) {};
    \node[draw = white] at (sum1.center) {$+$};
    
    \node [smallblock, fill=green!20, rounded corners, right = 1cm of sum1, minimum height = 0.6cm , minimum width = 1cm] (system) {$\SM$};
    \node [output, right = 0.5cm of system] (output) {};
    \node [input, below = 0.9cm of system] (midpoint) {};

    \draw[->] (estimator.east) -- node [above] {$\tilde{J}_{u}(t)$} (integrator.west);

    \draw[->] (integrator.east) -- (int_gain.west);
    \draw[->] (int_gain.east) -- node [above] {$\delta{u} (t)$} (sum1.west);
    \draw [->] (sum1) -- node [above] {$u (t)$} (system);
    \draw [->] (system) -- node [name=y, near end]{} node [very near end, above] {$J (t)$}(output);

    \draw[->] ([yshift = 0.85cm]sum1.north) -- node [xshift = 0.4cm, yshift = 0.4cm] {$d (t)$} (sum1.north);
    
    \draw [-] (y.west) |- (midpoint);
    \draw [->] (midpoint) -| ([xshift = -0.75cm]estimator.west) --  (estimator.west);
    
    \end{tikzpicture}
    }  
    \caption{Continuous-time extremum seeking control (ESC) of the continuous-time system $\SM.$
    ESC uses $J$ as an input and generates the control $u_k$ at each step $k$.
    Note that all components of $J$ are nonnegative.
    The objective of the controller is yield an input signal that minimizes the output of the continuous-time system, that is, yield $u(t)$ such that $\lim_{t \to \infty} J(t) = 0.$
    }
    \label{ESC_CT_blk_diag}
\end{figure}

Consider the system $\SM$ shown in first order extremum seeking scheme displayed in figure \ref{ESC_CT_blk_diag} to be
\begin{align}
    \dot x 
        &= f(x,u),
        \label{eq:x_dot_ESC}\\ 
    J(t) 
        &= \SQ(x),
        \label{eq:J_esc}
\end{align}
where $f: \mathbb{R}^{n} \times \mathbb{R}^{m} \to \mathbb{R}^{n}$ and $\SQ(x): \mathbb{R}^{n} \to \mathbb{R}^{p}$ are smooth enough, $x \in \mathbb{R}^n$ is the measured vector state, $u \in \mathbb{R}^m$ is the input vector and $J(t) \in \mathbb{R}^p$ is the output of the cost function $\SQ(x)$. Suppose that there exists $x^*$ such that $J^* = \SQ(x^*)$ is the extremum of the map $\SQ(.)$. Assume that both $x^*$ and $\SQ(.)$ are unknown. Thus, the main goal of extremum seeking control is to drive the states of the closed loop to $x^*$ without knowledge of $x^*$ or $\SQ(.)$

\subsection{SISO case}
Consider the case when $J(x) \in \mathbb{R}$, and $u \in \mathbb{R}$ are scalar. Next, consider the dither signal
\begin{align}
    d(t) = a~\sin(\omega_1~t),
\end{align}
where $a$ is the amplitude of the dither signal, and $\omega_1$ is the dither frequency.
Also, note that the gradient estimator used in this work is based on the averaging technique as proved in \cite{nevsic2009extremum} such that
\begin{align}
    \tilde{J}_u (t)
        &=
            \dfrac{2}{a}~J(t)\sin(\omega_1~t).
\end{align}
Finally, $\delta u(t)$, computed as the output of $\SI(t)$ shown in figure \ref{ESC_CT_blk_diag}, is obtained from the gradient descent scheme given by
\begin{align}
    \delta u(t) 
        &= 
            K_{\text{esc}}~\int_0^t \tilde{J}_u~dt.
\end{align}
In this scheme, $a,\omega_1$ and $K_{\text{esc}}$ are the tuning parameters.

\subsection{MISO case}
Now, consider the case when $J(x) \in \mathbb{R}$, and $u \in \mathbb{R}^m$. Define the vector of dither signals
$\SD(t)
    \isdef
        \matl
            d_1(t) &
            d_2(t) &
            \cdots &
            d_m(t)
        \matr^{\rm T}
$
given by
\begin{align}
   \SD(t) 
        &=
            a \matl
                \sin(\omega_1~t)\\
                \sin(\omega_2~t)\\
                \vdots\\
                \sin(\omega_m~t)
            \matr,
    \label{eq:MISO_dither}
\end{align}
where each $\omega_1,\omega_2,\cdots,\omega_m$ must be different. Note that although different amplitudes can be chosen for each dither signal, in the present work the same amplitude has been used for all, as shown in \eqref{eq:MISO_dither}. Define 
$\tilde{\SJ}_u 
    \isdef 
            \matl
                \tilde{J}_{u_1} &
                \tilde{J}_{u_2} &
                \cdots &
                \tilde{J}_{u_m}
            \matr^{\rm T}$.
Then, the gradient estimator based on the work done by \cite{ariyur2003real} is given by
\begin{align}
    \tilde{\SJ}_u
        &=
            \dfrac{2}{a}~J(t)
            \matl
               \sin(\omega_1~t)\\
               \sin(\omega_2~t)\\
               \vdots\\
               \sin(\omega_m~t)
            \matr.
\end{align}
Finally, Define
$\Delta u(t)
    \isdef
        \matl
           \delta u_1(t) &
           \delta u_2(t) &
           \cdots &
           \delta u_m(t) 
        \matr^{\rm T}$
given by the expression
\begin{align}
    \Delta u(t)
        &=
            K_{\text{esc}}\int_0^t\tilde{\SJ}_u~dt.
        \label{eq:MISO_GD_optimizer}
\end{align}
Note that although different $K_{\text{esc}}$ can be chosen for each component of $\tilde{\SJ}_u$, in the proposed scheme only one is considered, as shown in \eqref{eq:MISO_GD_optimizer}. Thus, $a$, $K_{\text{esc}}$ and 
$
    \matl
        \omega_1 &
        \omega_2 &
        \cdots &
        \omega_m
    \matr
$
are the tuning parameters.

\section{Overview of Retrospective Cost based Extremum Seeking Control} \label{sec:RCAC_ESC}

An overview of the RC/ESC algorithm is presented in this section.
Subsection \ref{RCAC:summary} presents a brief review of RCAC.
Subsection \ref{subsec:KF} describes an online gradient estimator based on the KF.
Subsection \ref{subsec:GDDRCAC} expands RCAC presented in Subsection \ref{RCAC:summary} to include gradient estimates obtained via the technique presented in Subsection \ref{subsec:KF}, resulting in RC/ESC.


\subsection{Review of Retrospective Cost Adaptive Control} \label{RCAC:summary}

RCAC is described in detail in \cite{rahmanCSM2017}. In this subsection we summarize the main elements of this method.
Consider the strictly proper, discrete-time, input-output controller  
\begin{equation}
    \delta u_k = \sum_{i=1}^{l_\rmc}P_{i,k} u_{k-i} + \sum_{i=1}^{l_\rmc}Q_{i,k}z_{k-i}, \label{IO_controller}
\end{equation}
where $k \ge 0$ is the time step, 
$\delta u_k \in \BBR^m$ is the RCAC controller output,
$u_k \in \BBR^m$ is the RC/ESC output and thus the control input,
$z_k \in \BBR^p$ is the measured performance variable, $l_\rmc$ is the controller-window length, and, for all $i\in \{1,\ldots,l_\rmc\},$  $P_{i,k} \in \BBR^{m \times m}$ and $Q_{i,k} \in \BBR^{m \times p}$ are the controller coefficient matrices.
In particular, $u_k$ results from adding a perturbation signal to $\delta u_k,$ as shown in Subsection \ref{subsec:GDDRCAC}.
The controller shown in \eqref{IO_controller} can be written as
\begin{equation}
    \delta u_k   =   \phi_k  \theta_k , \label{controller}
\end{equation}
where 
\begin{align}
	\phi_k &\isdef
\footnotesize		
\left[ \arraycolsep=3pt\def\arraystretch{0.9} \begin{array}{cccccc} 
			u_{k-1}^\rmT & \cdots & u_{k-l_\rmc}^\rmT & z_{k-1}^\rmT & \cdots & z_{k-l_\rmc}^\rmT
\end{array} \right]
		\otimes
		I_m
		\in \mathbb{R}^{m \times l_{\theta}},  \label{controller_phi} \\
	\theta_k &\isdef {\rm vec}
\left[ \arraycolsep=1.7pt\def\arraystretch{0.9} \begin{array}{cccccc} 
    P_{1,k} &\cdots &P_{l_\rmc,k} &Q_{1,k} &\cdots &Q_{l_\rmc,k}
\end{array} \right] \in \BBR^{l_\theta}, \label{controller_theta}
\end{align}
$l_\theta \isdef l_\rmc m (m + p),$ $\theta_k$ is the vector of controller coefficients, which are updated at each time step $k$,
and $\otimes$ is the Kroenecker product.

Next, define the retrospective cost variable
\begin{equation}
	\hat z_k (\hat \theta) \isdef z_k  - G_\rmf(\textbf{q})(\delta u_k - \phi_k \hat{\theta}), \label{zhat1}
\end{equation}
where $\hat{z}_k $ is the retrospective-cost variable,  
$G_{\rmf}$ is an $p \times m$ asymptotically stable, strictly proper transfer function, 
$\textbf{q}$ is the forward-shift operator,
and
$\hat{\theta} \in \BBR^{l_\theta}$ is the controller coefficient vector determined by  optimization below.
The rationale underlying \eqref{zhat1} is to replace the applied past control inputs with the re-optimized control input $\phi_k \hat{\theta},$ as mentioned in \cite{rahmanCSM2017} and \cite{islam2021data}.

In the present work, $G_{\rmf}$ is chosen to be a finite-impulse-response transfer function of window length $l_\rmf$ of the form
\begin{align}
	G_{\rmf}(\bfq) \isdef
			\sum_{i=1}^{l_\rmf} N_{i,k} \textbf{q}^{-i},
\end{align}
where  $N_{1,k},\ldots,N_{l_\rmf, k}\in\BBR^{p\times m}.$ 
We can thus rewrite \eqref{zhat1} as 
\begin{align}
    \hat z_k(\hat \theta) = z_k - N_k ( U_k - {\Phi}_k\hat{\theta} ),
\end{align}
where
\begin{align}{\Phi}_k \isdef 
    \footnotesize 
\left[ \arraycolsep=0.8pt\def\arraystretch{0.9} \begin{array}{c} 
        \phi_{k-1}       \\
         \vdots          \\
         \phi_{k-l_\rmf} \\
    \end{array}  \right] \in \BBR^{ l_\rmf m  \times l_\theta },\
    U_k \isdef 
    \footnotesize 
\left[ \arraycolsep=0.8pt\def\arraystretch{0.9} \begin{array}{c} 
        u_{k-1}       \\
        \vdots        \\
        u_{k-l_\rmf}  \\
    \end{array}   \right] \in \BBR^{ l_\rmf m   },
\end{align}
\begin{align}
    N_k \isdef \footnotesize \left[\  N_{1, k}  \ \cdots  \  N_{l_\rmf, k} \   \right] \in \BBR^{p \times l_\rmf m }.
\end{align}
In most applications, $N_k$ is constant and is determined by features of the system being controlled, as mentioned in \cite{rahmanCSM2017}. 
Other applications may require $N_k$ to be constructed and updated online using data, as mentioned in \cite{islam2021data}.
For the present work, the algorithm used to determine $N_k$ at each step $k$ is given in the next section.

\vspace{0.4em}

Using (\ref{zhat1}), we define the  cumulative cost function
\footnotesize
\begin{equation}
J_{\rmR, k}(\hat{\theta}) \isdef \sum\limits_{i=0}^{k} [  \hat z_i^{\rm T}(\hat \theta) \hat z_i(\hat \theta) +   (\phi_i \hat \theta)^{\rm T}  R_\rmu \phi_i \hat \theta ]  +   (\hat \theta - \theta_0 ) ^{\rm T}   P_0^{-1} (\hat \theta - \theta_0 ), \label{Jg}
\end{equation}
\normalsize
where $P_0 \in \BBR^{l_\theta \times l_\theta}$ is positive definite and $R_\rmu \in \BBR^{m \times m}$ is positive semidefinite. 
The following result uses recursive least squares (RLS), as mentioned in \cite{ljung:83} and \cite{islam2019recursive}, to minimize \eqref{Jg}, where, at each step $k,$  the minimizer of \eqref{Jg} provides the update $\theta_{k+1}$ of the controller coefficient vector $\theta_k$.

\textit{Proposition}. Let $P_0$ be positive definite, and $R_\rmu$ be positive semidefinite. Then, for all $k \ge 0$,  unique global minimizer $\theta_k$ of \eqref{Jg} is given by
\footnotesize
\begin{align}	
    P_{k}      = & \ P_{k-1}  - P_{k-1} 
    \left[ \arraycolsep=0.9pt\def\arraystretch{0.9} \begin{array}{c} 
                N_{k-1} {\Phi}_{k-1} \\
                \phi_{k-1} 
            \end{array}\right]^\rmT
    \hspace{-0.5em}\Gamma_{k-1}
    \left[ \arraycolsep=0.9pt\def\arraystretch{0.9} \begin{array}{c} 
                N_{k-1} {\Phi}_{k-1} \\
                \phi_{k-1} 
            \end{array}\right]
    P_{k-1} , \label{eq:pk_update}\\
    \theta_{k} = & \ \theta_{k-1} \nn \\
                & \hspace{-1em} - P_k             
\left[ \arraycolsep=0.9pt\def\arraystretch{0.9} \begin{array}{c} 
                N_{k-1} {\Phi}_{k-1} \\
                \phi_{k-1}
            \end{array}\right]^\rmT \hspace{-0.5em} \bar{R}             
\left[ \arraycolsep=0.9pt\def\arraystretch{0.9} \begin{array}{cc} 
                 z_{k-1} - N_{k-1} ( U_{k-1} - {\Phi}_{k-1} \theta_{k-1} ) \\
                \phi_{k-1} \theta_{k-1}
            \end{array}\right] \label{eq:theta_update},
\end{align}
\normalsize
where
\footnotesize
\begin{align}	
   \Gamma_{k-1} \isdef & \
                 \bar{R} - 
                    \bar{R} \left[ \arraycolsep=1.1pt\def\arraystretch{0.9} \begin{array}{c} 
                N_{k-1} {\Phi}_{k-1} \\
                \phi_{k-1} 
            \end{array}\right]
            \Xi_k^{-1}
            \left[ \arraycolsep=1.1pt\def\arraystretch{0.9} \begin{array}{c} 
                N_{k-1} {\Phi}_{k-1} \\
                \phi_{k-1} 
            \end{array}\right]^\rmT
            \bar{R} \nonumber \\
            &\in \BBR^{ (p + m) \times  (p + m) } ,\\
            \Xi_k \isdef & \  P_{k-1}^{-1}
            +
            \left[ \arraycolsep=1.1pt\def\arraystretch{0.9} \begin{array}{c} 
                N_{k-1} {\Phi}_{k-1} \\
                \phi_{k-1} 
            \end{array}\right]^\rmT
            \bar{R}
            \left[ \arraycolsep=1.1pt\def\arraystretch{0.9} \begin{array}{c} 
                N_{k-1} {\Phi}_{k-1} \\
                \phi_{k-1} 
            \end{array}\right] \in \BBR^{l_\theta \times l_\theta} , \\
    \bar{R} \isdef & \ {\rm diag}( I_p , R_\rmu ) \in \BBR^{ (p + m) \times  (p + m) }.
\end{align}
\normalsize

For all of the numerical simulations and physical experiments in this work, $\theta_k$ is initialized as $\theta_0=0_{l_\theta\times 1}$ to reflect the absence of additional prior modeling information.
The matrices $P_0$ and $R_\rmu$ have the form  $P_0 = p_0 I_{l_\theta}$ and $R_\rmu = r_\rmu I_m,$ where the positive scalar $p_0$ and nonnegative scalar $r_\rmu$ determine the rate of adaptation.

\subsection{RCAC-based PID}\label{subsec:RCAC_PID}
In the case where $m = p = 1,$ let $\delta u_k$ be given by
\begin{equation}
    \delta u_k=\kappa_{\rmp,k} z_{k-1}+\kappa_{\rmi,k} \zeta_{k-1}+\kappa_{\rmd,k}(z_{k-1} - z_{k-2}), \label{eq:pid}
\end{equation}
where $\kappa_{\rmp,k}$, $\kappa_{\rmi,k}$, and $\kappa_{\rmd,k}$ are time-varying PID gains 
%
%
%
and $\zeta_{k}$ is given by the integrator
\begin{align}
    \zeta_{k} \isdef  \sum_{j=0}^{k} z_j.
\end{align}
The control \eqref{eq:pid} can be expressed as
\begin{equation}
    \delta u_k=
      \phi_{k} \theta_{k}, \label{pid_controller}
\end{equation}
where
\begin{align}
    \phi_k \isdef &
        \matl
             z_{k-1} & 
            \zeta_{k-1} &
             z_{k-1} -  z_{k-2}  
        \matr \in \mathbb R^{1 \times 3}, \label{phithetadefns1} \\
    \theta_k  \isdef &
        \matl
         \kappa_{\rmp,k} &
         \kappa_{\rmi,k} &
         \kappa_{\rmd,k}
        \matr^\rmT \in \mathbb R^{3}.\label{phithetadefns2}
\end{align}
Note that the regressor $\phi_k$ is constructed from the past values of $z_k$ and $\zeta_k,$ and the controller coefficient vector $\theta_k$ contains the time-dependent proportional, integral, and derivative gains $\kappa_{\rmp,k},$ $\kappa_{\rmi,k},$ and $\kappa_{\rmd,k}$.
Furthermore, note that the adaptive digital PID control can be  specialized to adaptive digital PI, PD, ID, P, I, and D control by omitting the corresponding components of $\phi_k$ and $\theta_k.$
Then, RCAC-based PID (RCAC/PID) can be implemented by replacing \eqref{controller_phi} and \eqref{controller_theta} with \eqref{phithetadefns1} and \eqref{phithetadefns2}, respectively.


\subsection{Online gradient estimator using a Kalman Filter}\label{subsec:KF}
For all $k\ge0,$ let $J_k \isdef [J_{1,k} \ \ \cdots \ \ J_{l_J, k}]^\rmT \in \BBR^{l_J}$ be a cost function vector computed from system measurements, where, for all $i \in \{1, \ldots, l_J\},$ $J_{i,k} \ge 0$ is the $i$th component of $J_k,$ 
let $u_k$ be the  control input,
and let $\nabla J_k \isdef [\nabla J_{1,k} \ \ \cdots \ \ \nabla J_{l_J, k}]^\rmT \in \BBR^{l_J \times m}$ be the gradient of $J_k$ over $u_{k},$ where, for all $i \in \{1, \ldots, l_J\},$ the transpose of $\nabla J_{i,k} \in \BBR^{m}$ corresponds to the $i$th row of $\nabla J_k.$

Next, let $i \in \{1, \ldots, l_J\}.$
Consider the measurement model for $J_{i,k}$
\begin{equation} \label{eq:costModel}
    J_{i,k} = J_{\rmb, i} + \nabla J_{i,k}^\rmT u_{k},
\end{equation}
where $J_{\rmb, i} \in \BBR$ is a bias variable.  
Note that \eqref{eq:costModel} is an extension of (17) from \cite{gelbert2012}.
Furthermore, let $\nabla \hat{J}_{i, k} \in \BBR^{m}$ be an estimate of $\nabla J_{i,k},$ let $\hat{J}_{\rmb, i} \in \BBR$ be an estimate of $J_{\rmb, i},$ let $\hat{x}_{i,k} \isdef \begin{bmatrix} \nabla \hat{J}_{i,k}^\rmT & \hat{J}_{\rmb, i} \end{bmatrix}^\rmT  \in \BBR^{m + 1}$ be an estimate of $x_{i,k} \isdef \begin{bmatrix} \nabla J_{i,k}^\rmT & J_{\rmb, i} \end{bmatrix}^\rmT,$ and let $P_{i,k} \in \BBR^{(m + 1) \times (m + 1)}$ be the covariance of the estimate $\hat x_{i,k}$ of $x_{i,k}.$
Then, as indicated by \eqref{eq:costModel} and Section 3.1 of \cite{gelbert2012}, the estimate $\nabla \hat{J}_{i,k}$ can be obtained using a KF with state and measurement equations given by
\begin{align}
    x_{i, k+1} &= x_{i,k} + w_{i,k}, \label{eq:x_KF} \\
    y_{i,k} \isdef \begin{bmatrix} J_{i,k} \\ J_{i,k-k_1} \\ \vdots \\ J_{i,k - k_{m}} \end{bmatrix} &= \begin{bmatrix} u_{k-1}^\rmT & 1 \\ u_{k-1-k_1}^\rmT & 1 \\ \vdots & \vdots \\ u_{k-1-k_{m}}^\rmT & 1 \end{bmatrix} x_{i,k} + v_{i,k}, \label{eq:y_KF}
\end{align}
where $y_{i,k}\in \BBR^{m + 1}$ is the measurement vector and $w_{i,k}, v_{i,k} \in \BBR^{m + 1}$ are Gaussian random variables.
Hence, it follows from \eqref{eq:x_KF}, \eqref{eq:y_KF} that the estimate $\nabla \hat{J}_{i,k}$ is given by the recursive update of the KF, 
whose prediction and update equations are given, for $i\in \{1, \ldots, l_J\},$ by
\begin{align}
\hat{x}_{i,k} &= \hat{x}_{i,k-1} + K_{i, k-1} \left( G_{i, k-1} - H_{k-1} \hat{x}_{i,k-1} \right), \label{eq:KF_xhat} \\
P_{i,k} &= (I_{m + 1} - K_{i, k-1} H_{k-1})(P_{i, k-1} + Q_i), \label{eq:KF_P}\\
\nabla \hat{J}_{i, k} &= \left[ I_{m} \ \ 0_{m \times 1} \right] \hat{x}_{i, k}, \label{eq:KF_nablaJ}
\end{align}
where 
\footnotesize
\begin{align*}
G_{i,k-1} \isdef \begin{bmatrix} J_{i,k-1} \\ J_{i,k-1-k_1} \\ \vdots \\ J_{i,k - 1 - k_{m}} \end{bmatrix} \in \BBR^{m + 1},
\
H_{ k-1} \isdef & \ \begin{bmatrix} u_{k-1}^\rmT & 1 \\
u_{k-1-k_1}^\rmT & 1 \\
\vdots & \vdots \\
u_{k-1-k_{m}}^\rmT & 1\end{bmatrix}  \\
& \ \in \BBR^{(m + 1) \times (m + 1)},
\end{align*}
\begin{align*}
K_{i, k-1} \isdef & \ [(P_{i, k-1} + Q_i) H_{k-1}^\rmT] [H_{k-1} (P_{i,k-1} + Q_i) H_{k-1}^\rmT + R_i]^{-1} \\
& \ \in \BBR^{(m + 1) \times (m + 1)}, 
\end{align*}
\normalsize
$Q_i, R_i \in \BBR^{(m + 1)\times(m + 1)}$ are the constant weighting matrices, and $0< k_1 < \dots < k_{m}$ are indices.
The matrices $Q_i$ and $R_i$ determine the rate of estimation, and $k_1, \ldots, k_{m}$ are chosen to enhance the accuracy of the estimate $\hat{x}_{i, k}.$
Finally, the estimate $\nabla \hat{J}_k$ is given by
\begin{equation}
    \nabla \hat{J}_k \isdef \left[\nabla \hat{J}_{1,k} \ \ \cdots \ \ \nabla \hat{J}_{l_J, k}\right]^\rmT \in \BBR^{l_J \times m}.
\end{equation}

For all of the numerical simulations in the present work, $\hat{x}_{i, k}$ is initialized as $\hat{x}_{i, 0} = 0_{m\times 1}.$
The matrices $P_{i,0},$ $Q_i,$ and $R_i$ have the form  $P_{i,0} = p_{i,0} I_{m + 1},$ $Q_i = q_i I_{m + 1},$ and $R_i = r_i I_{m + 1},$ where the positive scalars $p_{i,0},$ $q_i,$ and $r_i$ determine the rate of estimation.


\subsection{Retrospective Cost based Extremum Seeking Control}\label{subsec:GDDRCAC}

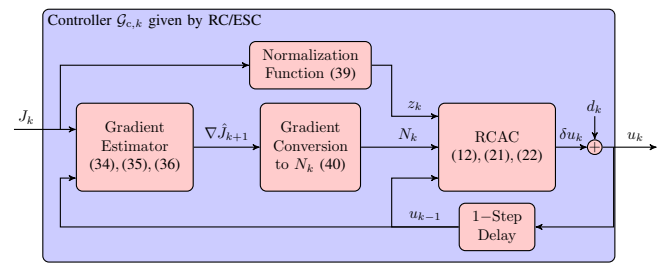
\begin{figure}[h!]
\centering
    \resizebox{\columnwidth}{!}{%
    	\begin{tikzpicture}[>={stealth'}, line width = 0.25mm]
		\node [smallblock,  rounded corners, minimum height = 5.75cm, minimum width = 13cm] (adap_controller) {};
		\node [input, name=ref, below left = -1.25cm and 8.75 cm of adap_controller.center]{};
		\node[below right] at (adap_controller.north west) {Controller $\SG_{\rmc, k}$ given by RC/ESC};
		\node [smallblock, fill=red!20,rounded corners, below right = 0.5cm and 3 cm of ref, minimum height = 2cm, minimum width = 2cm] (estimator) {$\begin{array}{c} {\rm Gradient} \\ {\rm Estimator} \\ \eqref{eq:KF_xhat},\eqref{eq:KF_P},\eqref{eq:KF_nablaJ} \end{array}$};
		\node [smallblock, fill=red!20,rounded corners, right = 1.5cm of estimator, minimum height = 2cm, minimum width = 2cm] (converter) {$\begin{array}{c} {\rm Gradient} \\ {\rm Conversion} \\ {\rm to} \ N_{k} \ \eqref{eq:grad_conv_eq} \end{array}$};
		\node [smallblock, fill=red!20,rounded corners, above = 0.3 cm of converter, minimum height = 0.75 cm, minimum width = 0.75cm] (normalization) {$\begin{array}{c} {\rm Normalization} \\ {\rm Function} \ \eqref{eq:norm_eq} \end{array}$};
		\node [smallblock, fill=red!20,rounded corners, right = 1.75cm of converter, minimum height = 2cm, minimum width = 2cm] (RCAC) {$\begin{array}{c} {\rm RCAC} \\ \eqref{controller},\eqref{eq:pk_update}, \eqref{eq:theta_update} \end{array}$};
        \node [smallblock, fill=red!20,rounded corners, below = 0.25cm of RCAC.south, minimum height = 1cm, minimum width = 1cm] (delay) {$\begin{array}{c} 1{\rm -Step} \\ {\rm Delay} \end{array}$};
		\node [sum, fill=red!20,right = 2em of RCAC.east] (sum1){};
		\node[draw = none] at (sum1.center) {$+$};
		\draw [->] ([xshift = -4em, yshift = 1.15em]estimator.west) -- node [above, near start, xshift = -0.2em] (J_k) {$J_{k}$} ([yshift = 1.15em]estimator.west);
		\draw [->] (estimator.east) -- node [above] (dJ_k) {$\nabla \hat{J}_{k+1}$} (converter.west);
		\draw [->] (converter.east) -- node [above, xshift = 0.55em] (Gf_k) {$N_{k}$} (RCAC.west);
		\draw [->] (RCAC.east) -- node [above] (u_rcac) {$\delta u_k$} (sum1.west);
		\draw [->] (sum1.east) -- node [above, near end, xshift = -0.5em] (u_req) {$u_k$}  ([xshift = 3.5em]sum1.east);
		\draw [->] ([xshift = -1em, yshift = 1.15em]estimator.west) |- (normalization.west);
		\draw [->] (normalization.east) -| ([xshift = -3em, yshift = 2em]RCAC.west) --  node [above] (zk) {$z_{k}$}  ([yshift = 2em]RCAC.west);
		\draw[<-] (sum1.north) -- node [near end, yshift = 1em] (u_exc) {$d_k$} ([yshift = 1.5 em]sum1.north);
        \draw[->] ([xshift = 0.25cm]sum1.east) |- (delay.east);
        \draw[->] (delay.west) -| node [near start, above] (ukm1) {$u_{k-1}$} ([xshift = -3em, yshift = -2em]RCAC.west) -- ([yshift = -2em]RCAC.west);
        \draw[->] (delay.west) -| ([xshift = -1em, yshift = -2em]estimator.west) -- ([yshift = -2em]estimator.west);
	\end{tikzpicture}
    }
\caption{RC/ESC block diagram.}
\label{quasi_blk_diag}
\end{figure}

As shown in Figure \ref{quasi_blk_diag}, RC/ESC includes RCAC described in Subsection \ref{RCAC:summary}, the KF gradient estimator described in Subsection \ref{subsec:KF}, a normalization function, and a gradient conversion function.
For RC/ESC, $l_J = p$ and $l_\rmf = m.$
RC/ESC operates on the cost-function vector $J_k \in \BBR^{p}$ and $u_{k-1} \in \BBR^{m}$ to produce the RC/ESC output vector $u_k \in \BBR^{m}.$
As mentioned in subsection \ref{subsec:KF}, for all $i \in \{1, \ldots, p\},$ $J_{i, k} > 0.$
The objective of RC/ESC is to minimize the components of $J_k$ by modulating $u_k,$ that is,
\begin{equation}
    \min_{(u_n)_{n=0}^{\infty}} \ \ \limsup_{k\to\infty} \sum_{i = 1}^{p} J_{i, k}.
\end{equation}

The performance variable $z_{k}$ used by RCAC is obtained by normalizing $J_{k}$ 
using
\begin{equation} \label{eq:norm_eq}
    z_{k} \isdef [I_{p} + \nu \ {\rm diag}\left( J_{k}\right)]^{-1} J_{k},
\end{equation}
where $\nu \in [0, \infty).$ 
%
%
Next, the gradient estimator block operates on $J_{k}$ and $u_{k}$ to produce $\nabla \hat{J}_{k+1}$ by using the KF-based estimator described in Subsection \ref{subsec:KF}.
The gradient conversion block yields $N_{k} = \left[ N_{1, k} \ \ \cdots \ \ N_{m, k} \right],$ such that,
for all $i \in \{1, \ldots, m\},$
\begin{equation} \label{eq:grad_conv_eq}
    N_{i, k} = \left[  \begin{matrix} & \overline{\nabla} \hat{J}_{1, i, k+1} & \\
    0_{p \times (i-1)} & \vdots &  0_{p \times (m - i)} \\
    & \overline{\nabla} \hat{J}_{p, i, k+1} & \end{matrix} \right],
\end{equation}
where, for all $j \in \{1, \ldots, p\},$
\begin{equation}
    \overline{\nabla} \hat{J}_{j, i, k+1} \isdef \begin{cases}
    \nabla \hat{J}_{j, i, k+1} / \norml \nabla \hat{J}_{j, k+1} \normr, & \norml \nabla \hat{J}_{j, k+1} \normr \ge \varepsilon \\
    \nabla \hat{J}_{j, i, k+1} / \varepsilon, & {\rm otherwise,}
    \end{cases}
\end{equation}
$\nabla \hat{J}_{j, i, k+1}$ is the $i$th component of $\nabla \hat{J}_{j, k+1},$ and $\varepsilon > 0.$
We fix $\varepsilon = 10^{-4}$ throughout the present work.
The RCAC block then uses $z_{k}, N_{k},$ and $u_{k-1}$ to produce $\delta u_{k} \in \BBR^{m}$ by using the operations described in Subsection \ref{RCAC:summary}.

Finally, define $u_k \isdef \delta u_k + d_k,$ where $d_k \in \BBR^{m}$ is a vanishing perturbation signal.
Note that, while \cite{gelbert2012}  uses only $\nabla \hat{J}_{k+1},$ RC/ESC uses $J_{k}$ and $\nabla \hat{J}_{k+1}$ in the form of $z_k$ and $N_k,$ respectively.

\section{Numerical Examples} \label{sec:numerical_examples}

In this section, RC/ESC is implemented in numerical simulations to illustrate its performance and compare it with the continuous-time ESC algorithms presented in Section \ref{sec:ESC}.
Example \ref{ex:static_siso} features a static optimization problem in a SISO system.
Example \ref{ex:static_miso} features a static optimization problem in a MISO system.
Example \ref{ex:dynamic_miso} features a dynamic optimization problem in a SISO system.
In Examples \ref{ex:static_siso} and \ref{ex:static_miso}, $T_\rms = 1$ s. In Example \ref{ex:dynamic_miso}, $T_\rms = 5$ s.
Tables \ref{example_tab_RC_ESC_hyper} and \ref{example_tab_ESC_hyper} show the RC/ESC and ESC hyperparameters, respectively, used in the numerical examples.
Note that in Example \ref{ex:static_siso} a RCAC/I controller is implemented, that is an RCAC-based integrator controller, as mentioned in Subsection \ref{subsec:RCAC_PID}. 
Examples \ref{ex:static_miso} and \ref{ex:dynamic_miso} implement the general case RCAC introduced in Subsection \ref{RCAC:summary}.

\begin{table}[h!]
\caption{RC/ESC hyperparameters in numerical examples}
\label{example_tab_RC_ESC_hyper}
\centering
\renewcommand{\arraystretch}{1.2}
\resizebox{\columnwidth}{!}{%
\begin{tabular}{|c|c|c|c|c|c|c|c|c|c|c|}
\cline{2-10}
 \multicolumn{1}{c|}{} & \multicolumn{4}{|c|}{\textbf{RCAC}} & \multicolumn{5}{|c|}{\textbf{KF}} & \multicolumn{1}{|c}{} \\
\hline
\textbf{Example} & \textbf{Type} & $\bm{l_{\rmc}}$ & $\bm{r_\rmu}$ & $\bm{p_0}$ & $\bm{p_{1,0}, p_{2,0}}$ & $\bm{q_1, q_2}$ & $\bm{r_1, r_2}$ & $\bm{k_1}$ & $\bm{k_2}$  & $\bm{\nu}$ \\
\hhline{|=|=|=|=|=|=|=|=|=|=|=|}
\ref{ex:static_siso} & RCAC/I & - & \multirow{2}*{0.05} & 0.9 & $10^{-3}$ & \multirow{2}*{0.1} & 10 & 3 & - & 0.9 \\
\cline{1-3}\cline{5-6}\cline{8-11}
\ref{ex:static_miso} & \multirow{2}*{RCAC}  & \multirow{2}*{5} &  & \multirow{2}*{0.1} & \multirow{2}*{$10^{-4}$} &  & 0.1 & \multirow{2}*{2} & \multirow{2}*{6} & \multirow{2}*{0.2} \\
\cline{1-1}\cline{4-4}\cline{7-8}
\ref{ex:dynamic_miso} &  &  & 0.01 &  &  & 0.01 & 1 &  & &   \\
\hline
\end{tabular}
}
\end{table}

\begin{table}[h!]
\caption{ESC hyperparameters in numerical examples}
\label{example_tab_ESC_hyper}
\centering
\renewcommand{\arraystretch}{1.2}
\resizebox{0.5\columnwidth}{!}{%
\begin{tabular}{|c|c|c|c|c|}
\hline
\textbf{Example} & $\bm{a}$ & $\bm{K_{\rm esc}}$ & $\bm{\omega_1}$ &  $\bm{\omega_2}$  \\
\hhline{|=|=|=|=|=|}
\ref{ex:static_siso} & 0.2 & \multirow{2}*{0.05} & 6 & - \\
\cline{1-2}\cline{4-5}
\ref{ex:static_miso} & 0.3 & & 30 & 50 \\
\hline
\ref{ex:dynamic_miso} & 0.2 & 5 & 3 & 5 \\
\hline
\end{tabular}
}
\end{table}

\begin{figure}[h!]
    \centering
    \vspace{0.5em}
    \resizebox{\columnwidth}{!}{%
    \includegraphics{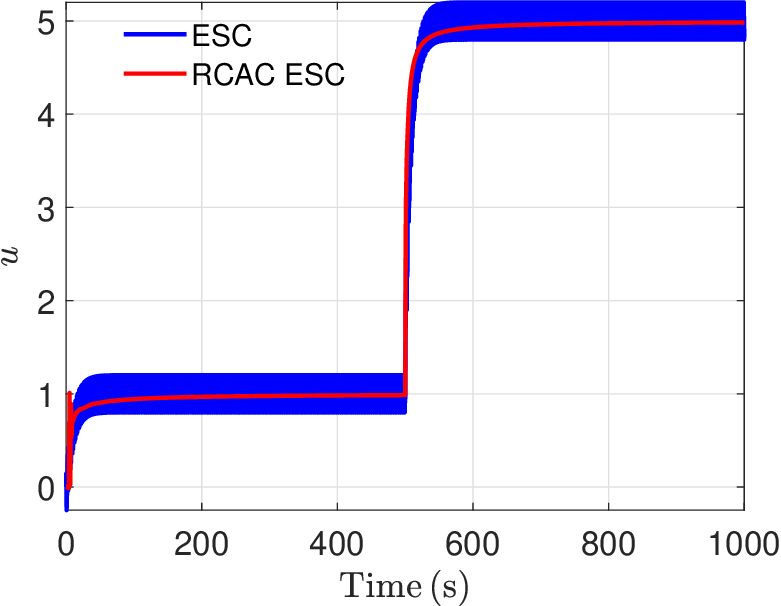}
    }
    \caption{Example \ref{ex:static_siso}:\textbf{SISO Static map}.  Controller output $u$ for the static map given by \eqref{SISO_static_example} with ESC and RC/ESC. Note that the ESC response is shown in blue and the RC/ESC response is shown in red.}
    \label{fig:ex1_u}
\end{figure}

\begin{figure}[h!]
    \centering
    \resizebox{\columnwidth}{!}{%
    \includegraphics{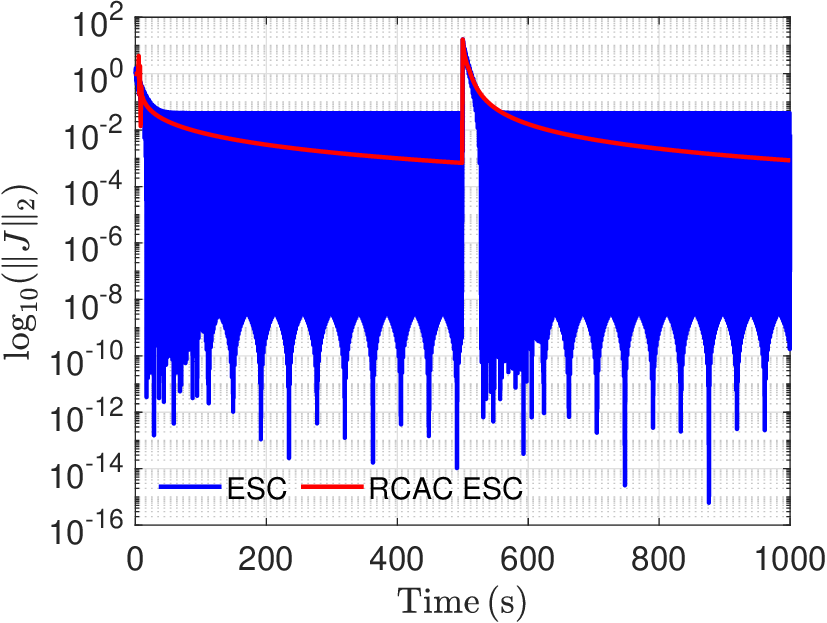}
    }
    \caption{Example \ref{ex:static_siso}:\textbf{SISO Static map}. Output error with respect to the optimal value $J=0$ in log scale with ESC and RC/ESC. Note that the error by ESC is shown in blue and the error by RC/ESC is shown in red.}
    \label{fig:ex1_logJ}
\end{figure}

\begin{figure}[h!]
    \centering
    \vspace{0.5em}
    \resizebox{\columnwidth}{!}{%
    \includegraphics{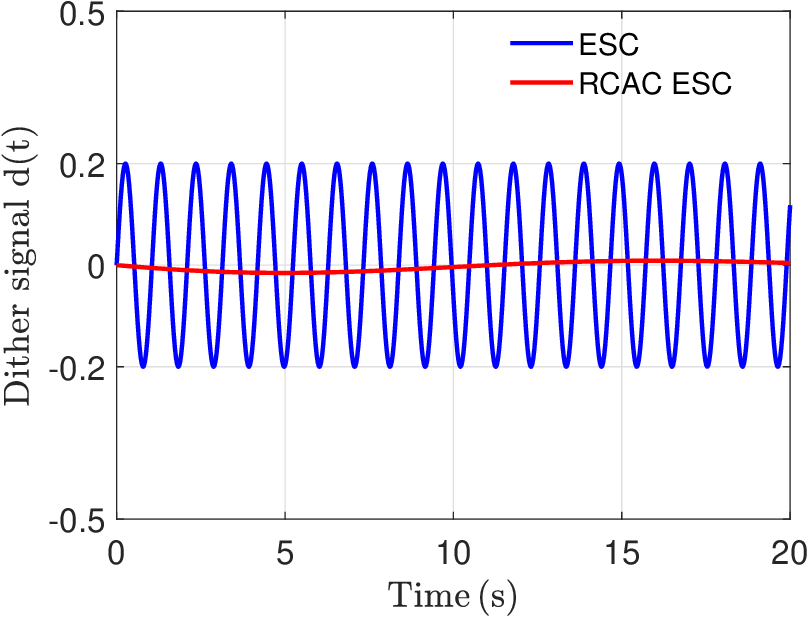}
    }
    \caption{Example \ref{ex:static_siso}:\textbf{SISO Static map}. Dither signal with ESC and RC/ESC. Note that the ESC dither signal is shown in blue and the RC/ESC dither signal is shown in red.}
    \label{fig:ex1_d}
\end{figure}

\begin{exam}\label{ex:static_siso}
{\it Static optimization in SISO system.}
Consider the static system
\begin{align}
    J(t) = (u(t) - r(t))^2
    \label{SISO_static_example}
\end{align}
where,
for $t \in [0,500], $ $r(t) = 1,$ and 
for $t > 500, $ $r(t) = 5.$
As mentioned in Section \ref{sec:control_statement}, the objective is to minimize $J.$
Furthermore, the dither signals are shown in \ref{fig:ex1_d}.
The results of the numerical simulations are shown in Figures \ref{fig:ex1_u} and \ref{fig:ex1_logJ}.
While the response of RC/ESC converges to the minimizer in a slow manner,
note that the response of RC/ESC doesn't oscillate around the minimizer since the dither is close to zero throughout the operation and eventually vanishes.
\hfill {\LARGE$\diamond$}
\end{exam}

\begin{figure}[h!]
    \centering
    \vspace{0.5em}
    \resizebox{\columnwidth}{!}{%
    \includegraphics{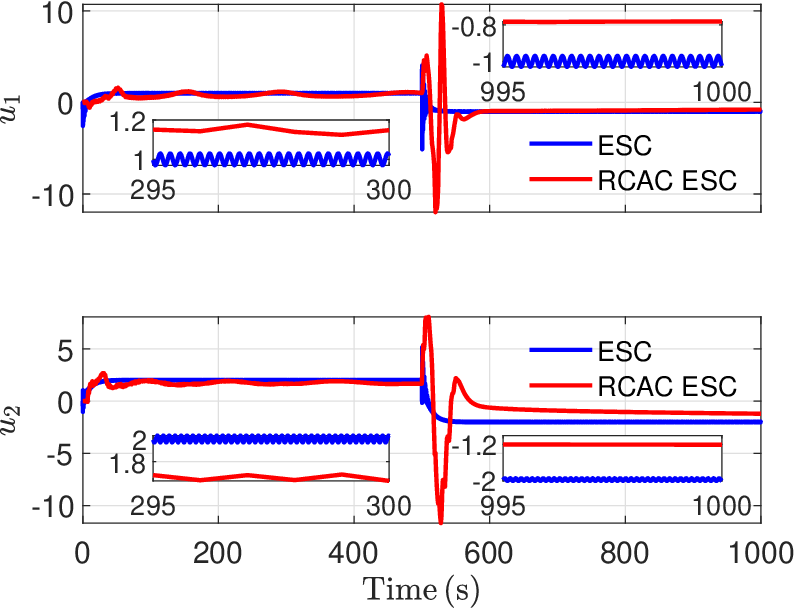}
    }
    \caption{Example \ref{ex:static_miso}:\textbf{MISO Static map}. Controller output components $u_{(1)}$ and $u_{(2)}$ for the static map given by \eqref{MISO_static_example} with ESC and RC/ESC. Note that the ESC response is shown in blue and the RC/ESC response is shown in red.}
    \label{fig:ex2_u}
\end{figure}

\begin{figure}[h!]
    \centering
    \vspace{0.5em}
    \resizebox{\columnwidth}{!}{%
    \includegraphics{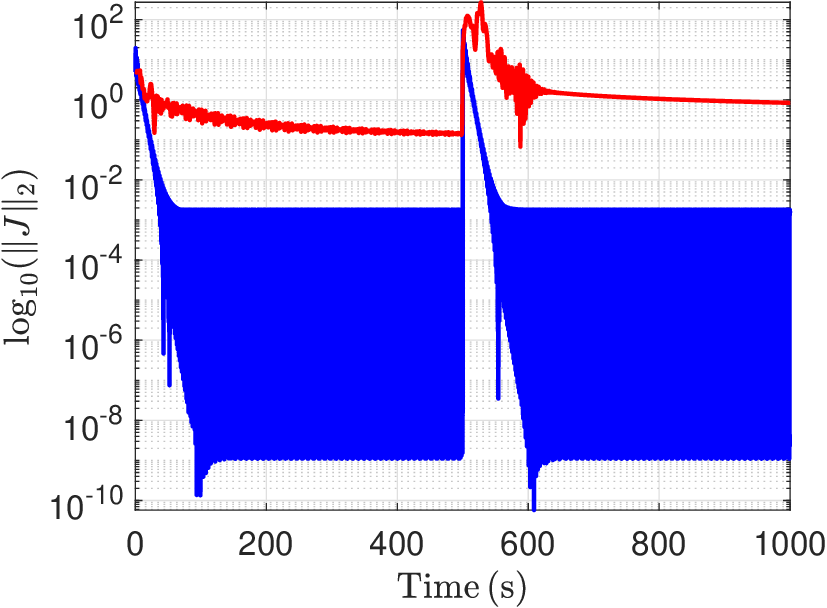}
    }
    \caption{Example \ref{ex:static_miso}:\textbf{MISO Static map}. Output error with respect to the optimal value $J=0$ in log scale with ESC and RC/ESC. Note that the error by ESC is shown in blue and the error by RC/ESC is shown in red.}
    \label{fig:ex2_logJ}
\end{figure}

\begin{figure}[h!]
    \centering
    \resizebox{\columnwidth}{!}{%
    \includegraphics{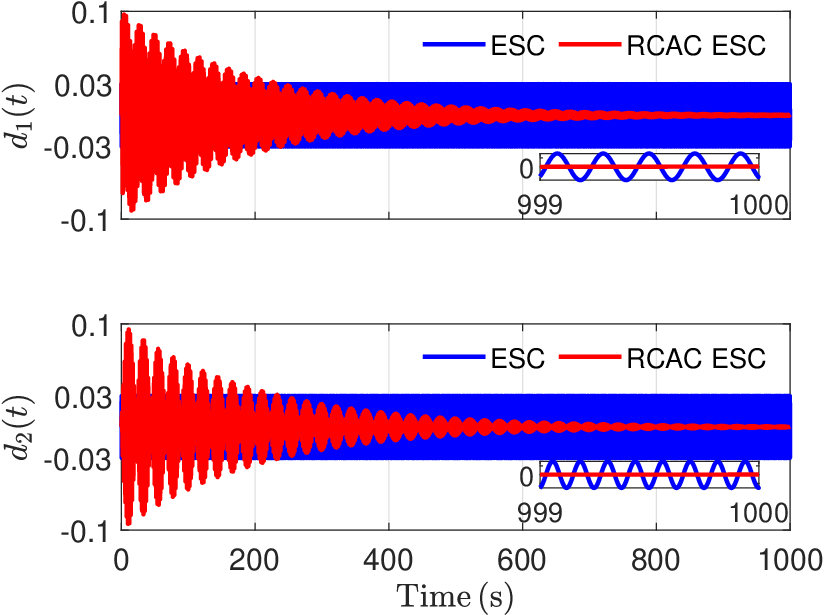}
    }
    \caption{Example \ref{ex:static_miso}:\textbf{MISO Static map}. Dither signals with ESC and RC/ESC. Note that the ESC dither signals are shown in blue and the RC/ESC dither signals are shown in red.}
    \label{fig:ex2_d}
\end{figure}

\begin{exam}\label{ex:static_miso}
{\it Static optimization in MISO system.}
 Consider the static system
\begin{align}
    J(t_ = (u_{(1)} (t) - r_{(1)} (t))^2 + (u_{(2)} (t) - r_{(2)} (t))^2
    \label{MISO_static_example}
\end{align}
where,
for $t \in [0,500], $ $r(t) = \matl
    1\\
    2
\matr,$ and 
for $t > 500, $ $r(t) = \matl
    -1\\
    -2
\matr.$
As mentioned in Section \ref{sec:control_statement}, the objective is to minimize $J.$
Furthermore, the dither signals are shown in \ref{fig:ex2_d}.
The results of the numerical simulations are shown in Figures \ref{fig:ex2_u} and \ref{fig:ex2_logJ}.
While the response of RC/ESC converges to the minimizer in a slow manner,
note that the response of RC/ESC doesn't oscillate around the minimizer since the dither quickly goes to zero throughout the operation and eventually vanishes.
\hfill {\LARGE$\diamond$}
\end{exam}

\begin{figure}[h!]
    \centering
    \vspace{0.5em}
    \resizebox{\columnwidth}{!}{%
    \includegraphics{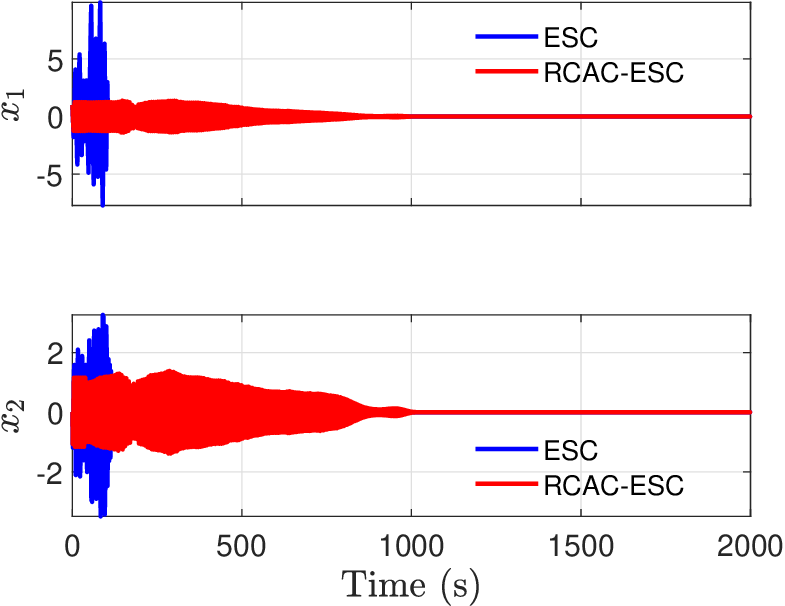}
    }
    \caption{Example \ref{ex:dynamic_miso}:\textbf{MISO Dynamic map}.
    Components of the Van Der Pol system state $x$ versus time with ESC and RC/ESC. Note that the ESC result is shown in blue and the RC/ESC result is shown in red.
    }
    \label{fig:ex3_xx}
\end{figure}

\begin{figure}[h!]
    \centering
    \vspace{0.5em}
    \resizebox{\columnwidth}{!}{%
    \includegraphics{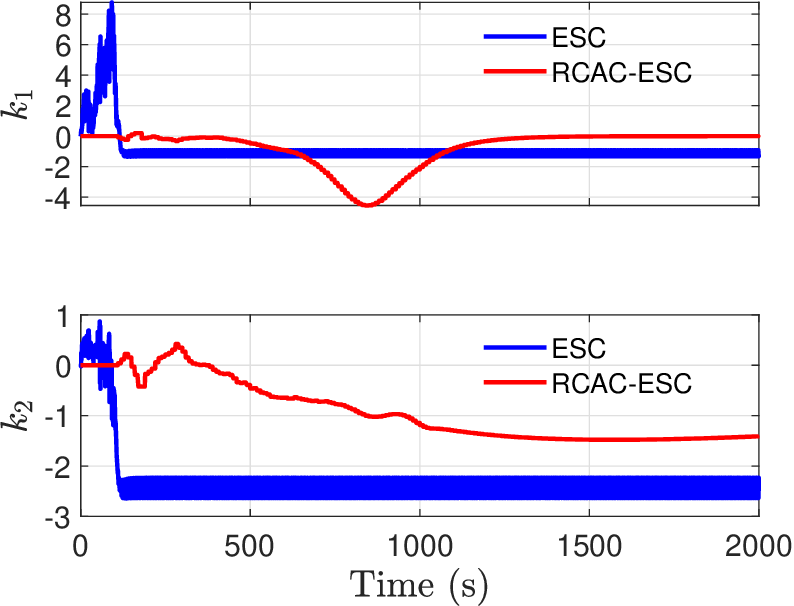}
    }
    \caption{Example \ref{ex:dynamic_miso}:\textbf{MISO Dynamic map}.
    Controller gains $K_1$ and $K_2$ versus time with ESC and RC/ESC. Note that the ESC response is shown in blue and the RC/ESC response is shown in red
    }
    \label{fig:ex3_u}
\end{figure}

\begin{figure}[h!]
    \centering
    \vspace{0.5em}
    \resizebox{\columnwidth}{!}{%
    \includegraphics{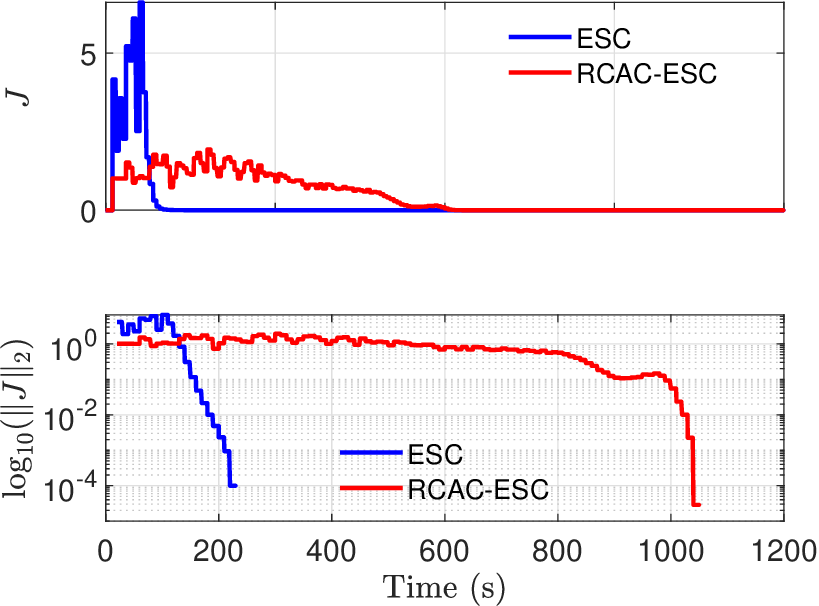}
    }
    \caption{Example \ref{ex:dynamic_miso}: \textbf{MISO Dynamic map}.
    Output error with respect to the optimal value $J=0$ in regular and log scale with ESC and RC/ESC. Note that the error by ESC is shown in blue and the error by RC/ESC is shown in red.
    }
    \label{fig:ex3_logJ}
\end{figure}

\begin{figure}[!ht]
    \centering
    \resizebox{\columnwidth}{!}{%
    \includegraphics{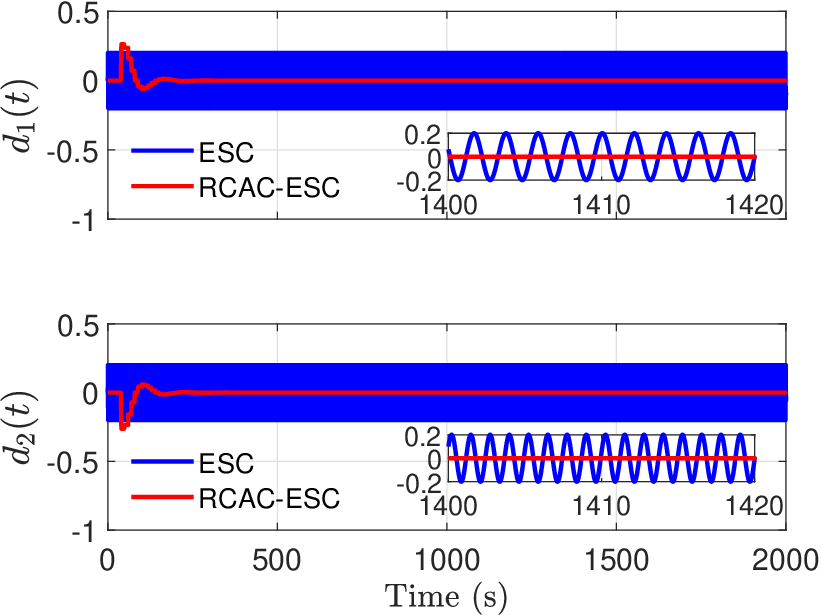}
    }
    \caption{Example \ref{ex:dynamic_miso}: \textbf{MISO Dynamic map}.
    Dither signals with ESC and RC/ESC. Note that the ESC dither signals are shown in blue and the RC/ESC dither signals are shown in red.
    }
    \label{fig:ex3_d}
\end{figure}

\begin{exam}\label{ex:dynamic_miso}
{\it Dynamic optimization in MISO system (control gain tuning for stabilization).}
\hfill \\
\\
Consider the Van Der Pol system
\begin{align}
    \ddot x + x + (\dot x^2 - 1)\dot x = u,
\end{align}
where $x_1 = x$ and $x_2 = \dot x$. Also, consider the full-state feedback controller structure given by
\begin{align}
    u 
        = 
            \matl 
                K_1 & K_2
            \matr
            \matl
                x_1\\
                x_2
            \matr.
\end{align}
Also, an amplitude detector scheme is considered using moving standard deviation for each state and adding them along the entire horizon. Thus, the cost function $J$ is the amplitude of the oscillations of the states and ESC and RCAC-ESC are used to find suitable values of $\matl K_1 & K_2 \matr$ in such a way that $J=0$ and thus, the system could be asymptotically stabilized.
As mentioned in Section \ref{sec:control_statement}, the objective is to minimize $J.$
Furthermore, the dither signals are shown in \ref{fig:ex3_d}.
The results of the numerical simulations are shown in Figures \ref{fig:ex3_xx}, \ref{fig:ex3_u} and \ref{fig:ex3_logJ}.
Both ESC and RC/ESC yield values of $K_1$ and $K_2$ that stabilize the response of the Van Der Pol system
While the response of RC/ESC converges to the minimizer in a slow manner,
note that the response of RC/ESC doesn't oscillate around the minimizer since the dither quickly goes to zero throughout the operation and eventually vanishes.
\hfill {\LARGE$\diamond$}
\end{exam}

\section{Conclusions}\label{sec:conclusions}

This paper introduced a retrospective cost-based ESC controller for online output minimization with a vanishing perturbation.
A KF is used to estimate the gradient of the system output at each step, which is then used to construct a target model that provides RCAC with a search direction to obtain a control input that minimizes the system output.
Numerical examples illustrate the performance of this technique and provide a comparison with a regular continuous-time ESC scheme.
Future work will focus on modifications for faster convergence and implementation in physical systems with time-varying minimizers.



\bibliographystyle{IEEEtran}
\bibliography{IEEEabrv,bib_paper}

\end{document}